\documentclass[12pt]{article}
\usepackage{graphicx}
\usepackage{amsmath,amsxtra,amssymb,latexsym, amscd,amsthm}
\usepackage[unicode]{hyperref}
\usepackage{cite}
\usepackage{enumitem}%gói lệnh giảm khoảng trống trong item Thêm vào trong document 2 dòng lệnh nữa mới được \setlist{noitemsep}\setlist{nolistsep}
\theoremstyle{plain}
\newtheorem{theorem}{Theorem}[section]
\newtheorem{definition}{Definition}[section]
\newtheorem{example}{Example}[section]
\newtheorem{corollary}{Corollary}[section]
\newtheorem{lemma}{Lemma}[section]

\newtheorem{remark}{Remark}[section]
\begin{document}
\setlist{noitemsep}
\setlist{nolistsep}
\title{On the Solution Existence  of  Nonconvex Quadratic Programming Problems in Hilbert Spaces}
\author{V. V. Dong\thanks{V. V. Dong, Phuc Yen College of Industry,  Vietnam,    vuvdong@gmail.com}        \and
        N. N. Tam \thanks{N. N. Tam, Hanoi Pedagogical Institute  2, Hanoi, Vietnam, nntam@hpu2.edu.vn }
} 
\date{05/10/2015}
\maketitle
\begin{quote}
\noindent {\bf Abstract.} {In this paper, we consider the  quadratic programming problems under  finitely many convex quadratic constraints in Hilbert spaces. 
By using the Legendre property of quadratic forms or the compactness of operators  in the presentations of quadratic forms, we establish some sufficient conditions for the solution existence of the considered problems. 
 As special cases, we obtain some existence solution results for   the quadratic programming problems under linear constraints in Hilbert spaces.}


\end{quote}
\section{Introduction}

Solution existence for quadratic programming problems (QP problems for brevity)  is an important and interesting question in optimization theory. This question has been studied extensively   both in the finite dimensional setting  and in the infinite dimensional setting.
% of both finite dimensional and infinite dimensional spaces.
 In 1956, Frank and Wolfe \cite {FW56}  proved the solution existence theorem for a finite dimensional QP problem. This result, called the Frank-Wolfe theorem, states that a quadratic function bounded below over a nonempty polyhedral convex set attains its infimum there. Since then, many extensions of the result have been obtained.   Alternative proofs and many improvings  of the Frank-Wolfe theorem were given by Eaves \cite{Eaves71}, Blum and Oettli  \cite {BO72}, Belousov \cite{Belousov77}, Luo and Zhang \cite{LuoZhang99}, Belousov and  Klatte \cite{BK02}. Bertsekas and Tseng  \cite{Bertsekas07} gave a series of existence results for the general quadratically constrained QP problems in the finite dimensional setting.  Semple \cite{Semple96},  Schochetman, Smith and Tsui  \cite{Schochetman97},  Bonnans and Shapiro \cite{BS00} studied  solution existence for  QP problems under generalized polyhedral constraint set  in  Hilbert spaces.  Borwein \cite{Borwein} investigated necessary and sufficient conditions for quadratic minimality.

The purpose of this paper is to obtain some new solution existence results for QP problem under finitely many convex constraints  in Hilbert spaces. Namely, in the following circumstances  the solution existence of  these problems will be established:

1) The quadratic form in the objective function is a Legendre form and  condition  \eqref{ConditionA} in section 3 is satisfied;

2)  The linear operators corresponding to quadratic forms describing the QP  problem are compact with closed range and the above-mentioned condition \eqref{ConditionA}  is satisfied.

The idea of using Legendre property of the quadratic form in objective function in proving solution existence of  quadratic programming problems with linear constraints in Hilbert space is due to Bonnans and Shapiro in \cite[Theorem 3.128]{BS00}. We would like to stress that the notion of Legendre form, which has origin in the Calculus of Variations is crucial for the just cited solution existence  theorem in \cite{BS00}. In section 3, we construct  an example to show that the conclusion of  that theorem fails if the assumption on the Legendre property of  the quadratic form is omitted.

By using an example from \cite{BK02} we can show that conditions \eqref{ConditionA} cannot be omitted if the given QP problem has more than one quadratic constraint. 

Note that the compact assumption in 2) guarantees  that the quadratic function is weakly continuous. 

The obtained results extend some preceding existence results for  nonconvex problems in infinite dimensional setting. As concerning the Hilbert space setting our results are completely new.

This paper is organized as follows. Some preliminaries are given  in Section 2. New existence results for the QP problem under finitely many quadratic constraints are derived in Section 3.

\section{Preliminaries}\label{sec:1}
In this section we recall some basic definitions and facts which we will use in the sequel. For details, we refer to \cite{BS00, Hestenes1, Ioffe, Hauser, Rudin}.

Let  $\cal H$ be a real Hilbert space with  inner product 
$\langle .  , . \rangle$ and the induced norm $\Vert .\Vert$.
A sequence $\{x_k\}$ in ${\cal H}$  is said to converge weakly to $x_0$, the notation  $x_k\rightharpoonup x$, if $\langle a, x_k\rangle\to \langle a, x_0\rangle$ for each $a$ in $\cal H$.  A sequence $\{x_k\}$ in ${\cal H}$  is said to converge strongly to $x_0$, the notation $x_k\to x$, if $\Vert x_k - x_0\Vert \to 0$.

A function $Q :  {\cal H} \to\mathbb R$  is said to be  quadratic  form  on $\cal H$  if  there exists a bilinear symmetric function  $B(. , .)$ on $\cal H$ such that  $Q(x)  =  B(x,  x)$.

In this paper, we will only consider  the continuous quadratic forms. By  Riesz Theorem (\cite[Theorem 2.34]{BS00}),  it admits the following representation:
\begin{align*}
Q(x) = \langle x, T x\rangle,
\end{align*}
where $T: {\cal H}\to {\cal H}$ is a continuous linear self-adjoint operator. 

 The operator $T:{\cal H}\to {\cal H}$ is said to be positive semidefinite (positive) if the quadratic from $Q(x) =\langle x, Tx\rangle$ is nonnegative (positive, respectively).

\begin{definition}{\rm (See, for instance, \cite[p.551]{Hestenes1})}
A quadratic form $Q(x)$ on the Hilbert space ${\cal H}$ is said to be a Legendre form if
\begin{itemize}
\item [{\rm (a)}]   it is weakly lower semicontinuous, and
\item [{\rm (b)}] $x_k\to x_0$ whenever $x_k\rightharpoonup x_0$ and $Q(x_k)\to Q(x_0)$.
\end{itemize}
\end{definition}

It is clear that in the case where ${\cal H}$ is of finite dimension, any quadratic form $Q(x)$ on ${\cal H}$ is a Legendre form.
\begin{example} Let $\ell^2$ denote the Hilbert space of all square summable real sequence. Define
$T: \ell_2\to \ell_2$ by  
$T x = (0, x_2, x_3,\ldots, x_n,\ldots)$, where $x=(x_1,x_2, x_3\ldots,x_n,\ldots)\in \ell^2$.
It is clear that $\langle x, T x\rangle = \frac{1}{2}(\Vert x\Vert^2  - x_1^2)$. By \cite[Proposition 3.79]{BS00},
 $\langle x, T x\rangle$ is a Legendre form.
\end{example}

\begin{example}\label{exam_notLeg}
	Let $L_2[0,1]$ denote the real Hilbert space of all square integrable function on $[0,1]$.  Let $T: L_2[0,1]\to L_2[0,1]$ be defined as follows: 
$$T x(t) = t x(t).$$
It is easy to check that $T$ is a  continuous linear self-adjoint and $\Vert T \Vert = 1$. The quadratic form associated with $T$ given by $Q(x) = \langle x, Tx \rangle = \int\limits_0^1 t x^2(t)dt$. Take
the sequence of functions $x_k: [0,1]\to \mathbb{R}$, defined for each positive integer $k$ by
\begin{align*}
x_k(t)=\begin{cases}
\sqrt{k}&\mbox{ if } \; 0\leqslant t\leqslant \frac{1}{k},\\
0& \mbox{ if }\; \frac{1}{k}\leqslant t\leqslant  1.
\end{cases}
\end{align*}
It is easy check that $x_k\in L_2[0, 1]$, $x_k$ converges weakly to zero and $x_k$ cannot converges strongly to zero. Furthermore,
\begin{align*}
Q(x_k)=\langle x_k, T x_k\rangle = \int\limits_0^1 t x_k^2(t)dt = k\int\limits_0^{\frac{1}{k}}tdt
=\frac{1}{2k}\rightarrow 0 = Q(0) \; (\mbox{ as } k\to\infty).
\end{align*}
Thus, $Q(x) = \langle x, T x\rangle$ is not a Legendre form.
\end{example}

A linear operator $A : {\cal H}\to {\cal H}$ is compact if it  transforms  every weakly convergent 
sequence into a strongly convergent sequence.

It well-known that of ${\cal H}$ is of finite dimension, then any continuous linear operator $T$ on ${\cal H}$  is compact.

It is easily seen that   the quadratic form $\langle x,Ix\rangle$ is a Legendre form while  $I$ is not a compact operator, the quadratic form  $\langle x,0x\rangle$ is not a Legendre form while  $0$ is  compact operator whenever $\cal H$ is of infinite dimension.

We will need the followings.
\begin{theorem}\label{S-lemma}(\cite{Hauser}, p.8)
Let ${\cal H}$ be a real Hilbert space, and let $g, h: {\cal H}\to \mathbb{R}$ be continuous quadratic functions defined on $\cal H$ by
\begin{align*}
g(x)=\frac{1}{2} \langle x, G x\rangle + \langle a, x\rangle +\gamma_1,\;\;\;
 h(x)= \frac{1}{2}\langle x,  Dx\rangle + \langle b, x\rangle +\gamma_2,
\end{align*}
and  {\it the Slater condition} is satisfied, that is there  exists
$x_0\in {\cal H}$ such that $h(x_0) < 0$.
Then $g(x)$ is bounded from below over $\{x\in {\cal H}\mid h(x) \leqslant 0\}$ if and if  the exists $\lambda \geqslant 0$ such that $g(x) + \lambda h(x) \geqslant 0$ for all $x\in {\cal H}$.
\end{theorem}

\section{Existence results}
\label{intro}
In what follows we are interested to obtain existence results for  quadratic programming problems of the following form
\begin{align}\label{eq:QP}
\begin{cases}
&\min f(x):=\frac{1}{2}\langle x, Tx\rangle  + \langle c, x\rangle\\
\mbox{s.t.}&x\in {\cal H}:\;   g_i(x):=\frac{1}{2}\langle x, T_i x\rangle + \langle c_i, x\rangle+\alpha_i\leqslant 0,\\
& \quad\qquad\quad\; i=1,2,\dots,m,
\end{cases}\tag{QP}
\end{align}
where $\cal H$ is a Hilbert space, $T:{\cal H}\to {\cal H}$ is continuous linear self-adjoint operator, $T_i$ is  positive semidefinite continuous linear self-adjoint operator on $\cal H$,  $c, c_i \in \cal H$, and $\alpha$, $\alpha_i$ are real numbers, $ i=1,2,\dots,m$.

If $T_i$ are zero operators for all $i=1,\dots,m$, then we say that \eqref{eq:QP} is a {\it quadratic programming problem under linear constraints} and denote it by ${\rm (QPL)}$.
 Note that if $T$ and $T_i$ are zero operators for all $i=1,\dots,m$, then \eqref{eq:QP} become a {\it linear programming problem }and denoted by ${\rm (LP)}$.

Let
\begin{equation*}
F=\{x\in {\cal H} \mid \; g_i(x)=\frac{1}{2}\langle x,T_ix \rangle+\langle c_i, x\rangle+\alpha_i \leqslant
0\ \, {\rm for\ all}\ \, i=1,...,m\}
\end{equation*}
  denote the constraint set of \eqref{eq:QP}.

For proving our  main results we need some lemmas.

\begin{lemma}\label{lema:1}
The constraint set $F$ of \eqref{eq:QP} is convex and weakly closed.
\end{lemma}

{\it Proof.} 
Since $T_i$ ($i=1,...,m$) is positive semidefinite, by Proposition 3.71 in \cite{BS00},   $g_i$  is convex. Hence  $F$  is closed and convex, and by  \cite[Theorem 2.23, p.24]{BS00}, $F$ is weakly closed.
\qed

\smallskip Recall that the {\it recession cone} of a nonempty closed convex set $X\subset {\cal H}$ is
defined (see \cite[p. 33]{BS00})  by
$0^+X =\{v\in {\cal H}\; \mid\; x+tv\in X \;\;  \forall   x\in X\;\;  \forall t\geqslant 0\}.$

\begin{lemma}\label{Recession Cone Formula_1}
If $F$ is nonempty, then
\begin{equation*}
0^+F=\{v\in {\cal H}\; \mid\; T_iv=0, \langle c_i, v\rangle\leqslant 0,\
\forall i=1,...,m\}.
 \end{equation*}
\end{lemma}
{\it Proof. } This proof is similar to the proof of Lemma 1.1 in \cite{KTY12}.
\qed
\begin{lemma}\label{levelset}
	Consider the problem \eqref{eq:QP}. Let $S_a=\{x\in {\cal F}\mid f(x)\leqslant a\}$, where $a\in \mathbb R$ is given. If $\langle x,Qx\rangle $ is weakly lower semicontinuous and  $S_a$ is nonempty then $S_a$ admits an element of minimal norm.
\end{lemma}
{Proof.}
	Suppose that $\langle x,Qx\rangle $ is weakly lower semicontinuous and  $S_a$ is nonempty. We show that $S_a$ has an element of minimal norm. Let $d= \inf \{\Vert y\Vert \mid y\in S_a\}$. Then, there exists a sequence $\{y^l\}$ in $S_a$ such that $\displaystyle d=\lim_{n\rightarrow\infty}\Vert y^l\Vert $. Since $\{y^l\}$ is bounded, it has a weakly convergent subsequence. Without loss of generality, we can assume that $y^l \rightharpoonup \bar y$ as $l\to\infty$. Since ${\cal F}$ is convex and closed,  by \cite[Theorem 2.23, p.24]{BS00}, it is weakly closed. Hence, $\bar y\in {\cal F}$. Since $\langle x, Qx\rangle$  is  weakly lower semicontinuous, one has
	\begin{align*}
	\langle \bar y, Q\bar y\rangle  &\leqslant \underset{l \to\infty}{\liminf}\langle y^l, Q y^l\rangle .
	\end{align*}
	From this it follows
	\begin{align*}
	f(\bar y)=\dfrac{1}{2}\langle \bar y, Q\bar y\rangle+\langle c, \bar y\rangle&  \leqslant  \underset{l \to\infty}{\liminf}\Big(\dfrac{1}{2}\langle y^l, Q y^l\rangle+ \langle c, y^l\rangle\Big)\leq a.
	\end{align*}
	Combining this with $\bar y\in {\cal F}$ we have that  $\bar y\in S_a$. This shows that $d=\Vert \bar y\Vert $ and $\bar y$ is a element of minimal norm of $S_a$. The proof of lemma is complete.
\qed

\begin{lemma}\label{lema:rec2}
 If $x^k\in F$ and $x^k\neq 0$ for all $k$ and $\Vert x^k\Vert^{-1}x^k$ weakly converges to $\bar v$ then $\bar v\in 0^+F.$
\end{lemma}
{\it Proof. }
Since $x^k\in F$, we have
\begin{align}\label{eqlema:rec2_1}
\frac{1}{2}\langle x^k,T_ix^k\rangle +\langle c_i,x^k\rangle+\alpha_i \leqslant 0\;\; i=1,...,m.
\end{align}
Multiplying both sides of  the inequalities in  \eqref{eqlema:rec2_1} by $\Vert x^k\Vert^{-2}$ and  letting $k\rightarrow \infty$, we obtain
 \begin{equation*}
  \underset{k\to\infty}{\liminf}\frac{1}{2}\left\langle \frac{x^k}{\Vert x^k\Vert}, T_i\frac{x^k}{\Vert x^k\Vert}\right\rangle \leqslant 0,\; i=1,...,m.
 \end{equation*}
 Since $T_i$ is positive semidefinite, by Proposition 3 in  \cite[p.~269]{Ioffe}, $\langle x,T_ix\rangle$ is  weakly lower semicontinuous.  Hence, 
 \begin{equation*}
\frac{1}{2}\langle \bar v, T_i\bar v\rangle \leq \underset{k\to\infty}{\liminf}
\frac{1}{2}\left\langle \frac{x^k}{\Vert x^k\Vert}, T_i\frac{x^k}{\Vert x^k\Vert}\right\rangle \leqslant 0\;\; i=1,...,m.
 \end{equation*}
 By  positive semidefiniteness of  $T_i$, from this we can deduce that
\begin{equation} \label{eqlema:rec2_2}
T_i\bar v=0\quad\, \forall i=1,...,m.
\end{equation}
As $\langle x^k,T_ix^k\rangle \geqslant 0$, from \eqref{eqlema:rec2_1} it follows
that
\begin{equation*}
\langle c_i, x^k \rangle+\alpha_i \leqslant 0,\quad \forall i=1,..., m,\ \, \forall k.
\end{equation*}
Multiplying the inequality $\langle c_i, x^k \rangle+\alpha_i \leqslant 0$ by
$\Vert x^k\Vert^{-1}$, letting $k\rightarrow \infty$, we get
\begin{equation}\label{eqlema:rec2_3}
\langle c_i,\bar v \rangle \leqslant 0\quad\, \forall i=1,...,m.
\end{equation}
Invoking Lemma \ref{Recession Cone Formula_1}, from  \eqref{eqlema:rec2_2},  
\eqref{eqlema:rec2_3} we have $\bar v\in {\rm Rec }(F)$.
\qed

\bigskip

 Consider the quadratic programming problem  \eqref{eq:QP}. Denote
\begin{equation*}
I=\{1,\dots, m\},\quad I_0=\{i\in I |\;  T_i=0\},\quad I_1=\{i\in I | \; T_i\not=0\}=I\setminus I_0 .
\end{equation*}
To prove our main results we need the following assumption:

\noindent{\bf Condition A}: If  $I_1\neq \emptyset$, then
\begin{align*}\label{ConditionA}
\left(v\in 0^+F,\langle v,Tv\rangle =0\right)\Rightarrow 
 \langle c_i,v\rangle =0\;\, \forall i\in I_1. \tag{A}
\end{align*}

Let us mention that in the finite dimensional setting, the condition \eqref{ConditionA} is  equivalent to the assumption  that all the asymptotic directions of $F$ are retractive local horizon directions  (see, \cite{Bertsekas07}). Furthermore, it is easily to see that 
 the \eqref{ConditionA} is satisfied if  one of the following conditions holds: 
\begin{itemize}
\item [\rm (1)] $T_i=0$ for all $i=1,\dots,m$
\item [\rm (2)] $\{v\in 0^+F\mid \langle v,Tv\rangle =0\}\subset \{0\}$
\item [\rm (3)] $c_i=0$ for all $i\in I_1$.
\end{itemize}

\begin{lemma}\label{lema:bd} Consider the problem \eqref{eq:QP}. Let $\{a^k\}$ be a   sequence of the real numbers,  and $S_k=\{x\in F \mid f(x)\leqslant a^{k}\}$. Suppose that
\begin{itemize}
\item [\rm (a)]  $\langle x, Tx\rangle$ is weakly lower semicontinuous;
\item [\rm (b)]  $x^k\in {\cal H}$ is a nonzero vector of minimal norm in $S_k$;
\item [\rm (c)]   $\frac{x^k}{\Vert x^k\Vert}$ weakly converges to nonzero $\bar v$ with $\langle \bar v, T \bar v\rangle = 0, \langle Tx^k+c, \bar v\rangle \geqslant 0$ for all $k$;
\item [\rm (d)] the condition \eqref{ConditionA} is satisfied.
\end{itemize}
Then, $\{x^k\}$  is bounded.
\end{lemma}
{\it Proof. } 
  By contrary, suppose that $\{x^k\}$ is unbounded.  Without loss of generality we may assume that $\Vert x^k\Vert \rightarrow \infty$ as $k\rightarrow \infty$, $\Vert x^k\Vert \neq 0$ for all $k$. Put $v^k:= \frac{x^k}{\Vert x^k\Vert}$, one has $\Vert v^k\Vert = 1$.  Since $\cal H$ is a Hilbert space, extracting if necessary a subsequence, we may suppose that $v^k $ itself weakly converges to some $\bar v$. 
From Lemma \ref{lema:rec2} it follows that 
\begin{equation}\label{eq:lema_bd1}
\bar v\in 0^+F.
\end{equation}

    We now show that there exists  $k_0$ such that $x^k-t\bar v\in F$ for all $k\geqslant k_0$ and for all $t>0$ small enough. To do this,  recall that $I=\{1,...,m\}$, $I_0=I\setminus I_1=\{i\mid T_i=0\}$ and $I_1=\{i\mid T_i\ne 0\}\subset I$.  Since $\{g_i(x^k)\}\leqslant 0$  
for all $i\in I$, we have the sequence $\{g_i(x^k)\}$ is bounded from above for each $i$. Therefore, there exists a subsequence $\{k^{\prime}\}$ of $\{k\}$ such that all the limits  $ \underset{k^{\prime}\to\infty}\lim g_i(x^{k^{\prime}})$ exist, $i=1,...,m$. Let us assume, without lost of generality, that $k^{\prime}=k$. Denote
\begin{align*}
I_{01}&=\{i\in I_0 \mid  \underset{k \to\infty}\lim g_i(x^{k})=0\}
=\{i\in I_0\mid \underset{k \to\infty}\lim (\langle c_i, x^{k}\rangle +\alpha_i)=0\}.\\
I_{02}&=\{i\in I_0 \mid  \underset{k \to\infty}\lim g_i(x^{k}) < 0\}
=\{i\in I_0\mid \underset{k \to\infty}\lim ( \langle c_i, x^{k}\rangle +\alpha_i) <0\}.
\end{align*}
Since $  \underset{k \to\infty}\lim (\langle c_i, x^{k}\rangle +\alpha_i )=0$  for all $i\in I_{01}$,  we can check that
\begin{equation*} 
\langle c_i,\bar v\rangle =0\;\; \forall i\in I_{01}.
 \end{equation*}
 This and condition \eqref{ConditionA} give
 \begin{equation}  \label{eq:lema_bd5}
\langle c_i,\bar v\rangle =0\;\; \forall i\in I_{01}\cup I_1.
 \end{equation}
 For $i\in I_{01}\cup I_1$, since $T_i=0$ for all $i\in I_0$, by  \eqref{eq:lema_bd1} and  \eqref{eq:lema_bd5}, we have
 \begin{align}  \label{eq:lema_bd6}
\begin{array}{ll}
 g_i(x^k-t\bar v)&= \frac{1}{2}\langle  x^k-t\bar v, T_i( x^k-t\bar v) \rangle+\langle c_i, x^k-t\bar v\rangle +\alpha_i\\
 &=\frac{1}{2}\langle  x^k, T_i x^k \rangle+\langle c_i, x^k\rangle +\alpha_i\\
 &=g_i(x^k) \leqslant 0\; \; \quad \quad \forall i\in  I_{01}\cup I_1, \forall\; t\geqslant 0.
\end{array}
  \end{align}
Since $ \underset{k \to\infty}\lim g_i(x^{k})= \underset{k \to\infty}\lim \langle c_i, x^{k}\rangle +\alpha_i <0$ for any  $i\in I_{02}$,  there exists  $\varepsilon >0$ such that
\begin{equation*}  
 \underset{k \to\infty}\lim g_i(x^{k})=\underset{k \to\infty}\lim (\langle c_i, x^{k}\rangle +\alpha_i) \leqslant - \varepsilon \; \quad  \forall i\in I_{02}.
 \end{equation*}
Let $i\in I_{02}$. Then, there exists $k_0$ such that
\begin{equation*}
g_i(x^{k})=\langle c_i, x^{k}\rangle +\alpha_i \leqslant -\frac{ \varepsilon}{2} \;\;\; \forall k\geqslant k_0.
 \end{equation*}
Fix $k\geqslant k_0$ and choose $\delta_{k,i}>0$ so that
$t\langle c_i, \bar v \rangle \geqslant -\frac{ \varepsilon}{2}$
for all $t\in (0,\delta_{k,i})$. Then, we have
 \begin{align}  \label{eq:lema_bd7}
\begin{array}{ll}
 g_i(x^k-t\bar v)&= \langle c_i, x^k-t\bar v\rangle +\alpha_i
 \leqslant \langle c_i, x^k\rangle +\alpha_i-t \langle c_i,\bar v\rangle\\
 &\leqslant   -\frac{ \varepsilon}{2}-t\langle c_i,\bar v\rangle
 \leqslant \; 0\quad \quad \forall i\in I_{02}.
\end{array}
  \end{align}
Let $\delta_k:=\min \{\delta_{k,i}\; |\; i\in I_{02}\}$. From \eqref{eq:lema_bd6}, 
\eqref{eq:lema_bd7}, it follows that
\begin{equation*} 
g_i(x^k-t\bar v)\leqslant 0\; \; \forall t\in (0,  \delta_{k})\quad \forall i=1,...,m.
\end{equation*}
This means that
\begin{equation} \label{eq:lema_bd8}
x^k-t\bar v\in F\quad \forall k\geqslant k_0\; \quad \forall  t\in (0,  \delta_{k}).
\end{equation}

We have, by assumption (iii),
\begin{align} \label{eq:lema_bd9}
\begin{array}{ll}
f(x^k-t\bar v)&=\frac{1}{2}\langle x^k-t\bar v, T(x^k-t\bar v)\rangle +\langle c, x^k-t\bar v\rangle\\
&= f(x^k)+\frac{t^2}{2}\langle \bar v,T\bar v\rangle-t\langle Tx^k+c,\bar v\rangle\\
&\leqslant  f(x^k).
\end{array}
\end{align}
Combining \eqref{eq:lema_bd8},  \eqref{eq:lema_bd9} we have
\begin{equation} \label{eq:lema_bd10}
x^k-t\bar v \in S_k \quad \forall k\geqslant k_0,\;  t\in (0,\delta_k).
\end{equation}
Since $\langle \bar v,\bar v\rangle =\Vert \bar v\Vert^2 >0$ and 
$\frac{x^k}{\Vert x^k\Vert}=v^k \rightharpoonup \bar v$, there exists $k_1\geqslant k_0$ such that
\begin{equation*}
\langle x^k,\bar v\rangle >0\quad \forall k\geqslant k_1.
\end{equation*}
Therefore, there exists $\gamma >0$ such that
\begin{align} \label{eq:lema_bd11}
\begin{array}{ll}
\Vert x^k-t\bar v\Vert^2&=\Vert x^k\Vert ^2-2t\langle x^k,\bar v\rangle +t^2\Vert\bar v\Vert ^2
<\Vert x^k\Vert^2\quad \forall t\in (0,\gamma).
\end{array}
\end{align}
  Let $\delta :=\min \{\delta_{k}, \gamma\}$. Then, by \eqref{eq:lema_bd10}, \eqref{eq:lema_bd11}, we have
  \begin{align*}
  &x^k-t\bar v\in S_k\quad \text{and \, }
   \Vert x^k-t\bar v\Vert <\Vert x^k\Vert  \quad \forall k\geqslant k_1, \quad \forall t\in (0,\delta).
    \end{align*}
   This contradics the fact that $x^k$ is the element of minimal norm in $S_k$.  Therefore, we conclude that $\{x^k\}$ must be bounded.
   
\qed

%\subsection{Existence theorem using the Legendre property}

 We now will give some existence results by using the Legendre property of the quadratic form in the objective function and the condition \eqref{ConditionA}.
\begin{theorem}\label{theo:exist1}{\rm (Frank-Wolfe-type theorem 1)}
 Consider the problem \eqref{eq:QP},
where $\langle x, Tx\rangle$ is a Legendre form.  Suppose that  $f(x)$ is bounded from below over nonempty $F$ and the condition \eqref{ConditionA} is satisfied.
 Then,  problem \eqref{eq:QP} has a solution.
\end{theorem}
{\it Proof. } Let $f^*=\inf \{f(x)\mid x\in F\}$. For each positive integer $k$, put $S_k=\{x\in F \mid f(x)\leqslant f^*+\frac{1}{k}\}$.
By $f$ is continuous and $f^* >-\infty$, $S_k$ is nonempty and closed. Since $\langle x,Tx\rangle$ is a Legendre form, it is weakly lower semicontinuous, by Lemma \ref{levelset}, $S_k$ admits an element of minimal norm, say  $x^k\in S_k$. Then, we have
\begin{align}
\label{eq:theo1}\; f(x_k)&=\frac{1}{2}\langle  x^k, Tx^k\rangle+\langle c, x^k\rangle   \leqslant f^*+\frac{1}{k},\\
 \label{eq:theo2} \; g_i (x_k)&=\frac{1}{2}\langle  x^k, T_i x^k\rangle+\langle c, x^k\rangle +\alpha_i  \leqslant 0\; \; i=1,\dots, m.
\end{align}
  Consider sequence $\{x^k\}$. We now prove that $\{x^k\}$ is bounded. On the contrary, suppose that $\{x^k\}$ is unbounded.  Without loss of generality we may assume that $\Vert x^k\Vert \rightarrow \infty$ as $k\rightarrow \infty$, $\Vert x^k\Vert \neq 0$ for all $k$. Put $v^k:= \frac{x^k}{\Vert x^k\Vert}$, one has $\Vert v^k\Vert = 1$, and then there exists a subsequence of $v^k$ weakly converges  to $\bar v$. Without loss of generality we may assume that $v^k \rightharpoonup \bar v\; \; \text {as}\; \; k\rightarrow \infty$.
From Lemma \ref{lema:rec2} it follows that $\bar v\in 0^+F$.

We are going to show that
\begin{align}
  \label{eq:theo3a}&\langle \bar v, T\bar v\rangle = 0\quad \text{and}\quad  \langle Tx^k+c,\bar v\rangle \geqslant 0\; \; \mbox{for each}\; k.
\end{align}
Multiplying the inequality in \eqref{eq:theo1}  by $\Vert x^k\Vert^{-2}$ and letting $k\to\infty$, by weakly lower semicontinuity of $\langle x,Tx\rangle$, we can deduce that
 \begin{equation}\label{eq:theo3b}
 \langle \bar v, T\bar v\rangle \leqslant \underset{k\to\infty}{\liminf}\frac{1}{2}\langle v^k, Tv^k\rangle \leqslant 0.
 \end{equation}
 Since $\bar v\in 0^+F$, by \eqref{eq:theo3b}, we must have \eqref{eq:theo3a}.
  Otherwise, for fix $k$ and for all $t>0$, $x^k+t \bar v\in F$, by \eqref{eq:theo3b} we have
 \begin{equation*}
f(x^k+t \bar v)= f(x^k) +\frac{t^2}{2} \langle \bar v,T\bar v\rangle +t\langle Tx^k+c, \bar v\rangle \to -\infty \; \mbox{as}\; t\to +\infty,
 \end{equation*}
  contradics the fact that  $f$ is  bounded from below over $F$.

We now show that $\Vert \bar v\Vert \neq 0$. Multiplying both sides of the inequality
 in \eqref{eq:theo1}  by $\Vert x^k\Vert^{-2}$ and letting $k\to\infty$, one has
 \begin{equation}\label{eq:theo3c}
  \underset{k\to\infty}{\limsup}\frac{1}{2}\langle v^k, Tv^k\rangle \leqslant 0.
   \end{equation}
Combining   \eqref{eq:theo3b},\eqref{eq:theo3c} and  $\langle \bar v,T\bar v\rangle=0$  we can conclude that
   \begin{equation}\label{eq:theo4}
    \underset{k\to\infty}{\lim}\langle v^k, Tv^k\rangle =\langle \bar v, T\bar v\rangle .
    \end{equation}
Since $\langle x,Tx\rangle$ is a Legendre form, by  $v^k \rightharpoonup \bar v$ and \eqref{eq:theo4}, we have $ v^k\rightarrow \bar v$. Since $\Vert v^k\Vert = 1$ for all $k$, we conclude that
   $\Vert \bar v\Vert =1$, so  $\bar v\neq 0$. 

 As $\langle x,Tx\rangle$ is a Legendre form, it is weakly lower semicontinuous. From the above it follows that all the conditions in  Lemma \ref{lema:bd} are satisfied. Hence by Lemma \ref{lema:bd}, $\{x^k\}$  is bounded. Since  $\{x^k\}$  is bounded, it has a weakly convergent subsequence. Without loss of generality, we can assume that $x^k\rightharpoonup \overline{x}$ as $k\to\infty$. 
Since $x_k\in F$ and $ F$ is weakly closed (see Lemma \ref{lema:1}), we have $ \overline{x}\in F $.
Since $\langle x, Tx\rangle $ is a Legendre form, it is  weakly lower semicontinuous,  one has
\begin{align*}
\frac{1}{2}\langle \bar x, T\bar x\rangle&\leqslant \underset{k\to\infty}{\liminf}\frac{1}{2}\langle x^k, Tx^k\rangle.
\end{align*}
Hence,  by  \eqref{eq:theo1},
\begin{align*}
f(\bar x)=\frac{1}{2}\langle \bar x, T\bar x\rangle+\langle c, \bar x\rangle&  \leqslant \underset{k\to\infty}{\liminf}\Big(\frac{1}{2}\langle x^k, T x^k\rangle+ \langle c, x^k\rangle\Big)\\
&\leqslant \underset{k\to\infty}{\liminf}(f^*+\frac{1}{k})
= f^*.
\end{align*}
It follows that  $\bar x$ is a solution of \eqref{eq:QP}. The proof is complete.

\qed

We now give an example  to show an application of  Theorem \ref{theo:exist1}.
\begin{example}\label{exam}
Let  ${\cal H}= \ell^2$,
$T x = (-x_1, 0, x_3,0, x_5,\ldots, x_n,\ldots)$, $T_1 x = ( x_1, 0, 0,\ldots)$, $T_2 x = ( 0, 0, 0, x_3, x_4, \ldots)$ for $x=(x_1,x_2,x_3 ,\ldots, x_n,\ldots)\in \ell^2$, and let
$c = (0,1,1,0\ldots), \, c_1 = (0,-1,-1,0,\ldots), \, c_2 =(0,0,\ldots) \in \ell^2$ , $\alpha_1 = 1, \alpha_2 = 0$. Then,
\eqref{eq:QP} becomes
\begin{align}\label{Prob_exam4}
\begin{cases}
&\mbox{min } f(x) = \frac{1}{2} (-x_1^2+0x_2^2+x_3^2+0x_4^2+x_5^2+\ldots) + x_2+x_3\\
&\mbox{subject to }   x\in F,
\end{cases}
\end{align}
where $F = \{x = (x_1,x_2,\ldots)\in \ell_2\mid  \frac{1}{2}x_1^2 -(x_2+x_3)+ 1\leqslant 0, x_3^2\leqslant 0\}$. 

It is clear that $F\neq\emptyset$ and  $-\frac{1}{2}x_1^2 +x_2+x_3\geqslant 1$ for any $x\in F$. 

Since 
\begin{align*}
f(x)&=\frac{1}{2} (0x_2^2+x_3^2+0x_4^2+x_5^2+\ldots) - \frac{1}{2} x_1^2+ x_2+x_3\\ 
&\geqslant 1 \;\;\; \mbox{for any \; } x\in F.
\end{align*}
Thus $f$ is bounded from below over $F$.
It is easy to check that 
$$\langle x, T x\rangle= \frac{1}{2}\Vert x \Vert^2 -\frac{1}{2}(2x_1^2+x_2^2+x_4^2).$$
 By \cite[Proposition 3.79]{BS00}, $\langle x, T x\rangle$  is a Legendre form. 
 According to Lemma \ref{Recession Cone Formula_1}, we have
\begin{align*}
0^+F = \{v=(0, v_2, 0, v_4, v_5,\ldots)\in\ell^2 \mid v_2\geqslant 0\}. 
\end{align*}
From this it follows that
$\{v\in 0^+F\mid \langle v, T v\rangle = 0\}=\{v=(0, 0, 0, v_4, 0,\ldots)\in\ell^2\}$.
Note that if $v\in \{v\in 0^+F\mid \langle v, T v\rangle = 0\}$, then $\langle c_1, v\rangle  = 0$. Thus the necessary condition \eqref{ConditionA} is satisfied. Take  $x^* = (1, 2, 0,\ldots)\in F$, we have   $f (x^* ) =1$. Since $f(x) \geqslant 1$ for all $x\in F$, $x^*$ is a  solution of \eqref{Prob_exam4}.
\end{example}

Let us mention some important consequences of the theorem.

\begin{corollary}\label{cor:exist1} (cf. \cite[Theorem 3.128]{BS00}) Consider the quadratic programming problem under linear constraints ${\rm (QPL)}$ {\rm(}i.e. \eqref{eq:QP} with $T_i=0$ for all $i=1,...,m${\rm)}, where $\langle x, Tx\rangle$ is a Legendre form. Suppose that $f(x)$ is bounded from below over nonempty $F$.
 Then,  problem ${\rm (QPL)}$  has a solution.
\end{corollary}
{\it Proof. } Since $T_i=0$ for all $i=1,...,m$, the set $I_1=\emptyset$. Hence the condition \eqref{ConditionA} is automatically satisfied and the corollary follows.
\qed
\begin{corollary}\label{cor:exist2} Consider the problem \eqref{eq:QP}, where $\langle x, Tx\rangle$ is a Legendre form. Suppose that  $c_i=0$ for all $i\in I_1$ and $f(x)$ is bounded from below over nonempty $F$.
 Then,  problem \eqref{eq:QP} has a solution.
\end{corollary}
{\it Proof. }
  Since $c_i=0$ for all $i\in I_1$, $\langle c_i,v\rangle =0$  for all $i\in I_1$. Hence the condition \eqref{ConditionA} is satisfied and the corollary follows.
\qed
\begin{corollary}\label{cor:exist3} Consider the problem \eqref{eq:QP}, where $\langle x, Tx\rangle$ is a Legendre form. Suppose that $\{v\in 0^+F\mid \langle v,Tv\rangle =0\}\subset \{0\}$  and $f(x)$ is bounded from below over nonempty $F$.
 Then,  problem \eqref{eq:QP} has a solution.
\end{corollary}
{\it Proof. }
  Since  $\{v\in 0^+F\mid \langle v,Tv\rangle =0\}\subset \{0\}$, $\langle c_i,v\rangle =0$  for all $i=1,...,m$ and the condition \eqref{ConditionA} is satisfied. The conclusion follows.
\qed
\begin{corollary}\label{cor:exist4}
Let $\langle x, Tx\rangle$  be a Legendre form on ${\cal H}$. Suppose that quadratic function  $f(x)=\frac{1}{2}\langle x, Tx\rangle+\langle c, x\rangle$ is bounded from below over Hilbert space ${\cal H}$. Then, there exists an $x^*\in {\cal H}$ such that $f(x^*)\leqslant f(x)$ for all $x\in {\cal H}$.
\end{corollary}
{\it Proof.} Consider   \eqref{eq:QP}  with $T_i=0$, $c_i=0$ and $\alpha_i=0$ for all $i=1,...,m$. Then, $F={\cal H}$ and it is clear that the condition \eqref{ConditionA} is satisfied. The conclusion follows.

\qed

The following example is constructed to show that the assumption on the Legendre property of  the quadratic form cannot be dropped from the assumptions of Theorem \ref{theo:exist1}.

\begin{example}\label{ex_nexist}
	Let $L_2[0,1]$ denote the real Hilbert space of all square integrable functions
	on $[0; 1]$ with the scalar product
	$$\langle x,y\rangle=\int\limits_0^1x(t)y(t)dt\quad \forall x,y\in L_2[0,1].$$
	
	Consider the programming problem (QP):
	
	\begin{equation}\label{QP}
	\begin{cases}
	\min f(x)=\dfrac{1}{2}\langle x, Tx\rangle\\
	\text {subject to} \quad x\in L_2[0, 1]:  \langle -c_1(t), x(t)\rangle +1\leqslant 0 
	\end{cases}
	\end{equation}
	where  $T: L_2[0,1]\to L_2[0,1]$ be defined by 
	$T x(t) = t x(t)$, and $c_1:[0,1]\to \mathbb{R}$, defined by
	\begin{equation*}
	c_1(t)=
	\begin{cases}
	\dfrac{1}{\sqrt{t}}\quad \text{if}\quad 0< t\leqslant 1\\
	0\quad \text{if}\quad t=0.
	\end{cases}
	\end{equation*}
It is easy to check that $T$ is a  continuous linear self-adjoint and $c_1(t)\in L_2[0,1]$.

Let $Q(x) = \langle x, Tx \rangle = \int\limits_0^1 t x^2(t)dt$. Then,  $Q(x)$ is not a Legendre form (see Example \ref{exam_notLeg})

Let
	$$F=\{x(t)\in L_2[0,1]\mid \; \langle -c_1(t), x(t)\rangle +1\leqslant 0 \}.$$

The set $F$ is nonempty. Indeed, consider the sequence of functions $x_n: [0,1]\to \mathbb{R}$, defined for each positive integer $n$ by

\begin{equation}
x_n(t)=
\begin{cases}
\sqrt {n} \;\quad \text{if}\;\; \dfrac{1}{n^2}\leqslant t \leqslant \dfrac{1}{n}\\
0  \;\;\quad \text{if}\;\; 0\leqslant t<\dfrac{1}{n^2} \; \; \text{or}\;\; \dfrac{1}{n}< t \leqslant 1.
\end{cases}
\end{equation}
It is easy to  see that $x_n\in L_2[0,1]$ and we have
\begin{equation}
\langle -c_1, x_n\rangle =-\int \limits_0^1c_1(t)x_n(t)dt=-\int \limits_{\dfrac {1}{n^2}}^{\dfrac{1}{n}}\dfrac{1}{\sqrt{t}}\sqrt{n}dt=-2\sqrt{n}(\dfrac{1}{\sqrt{n}}-\dfrac{1}{n})=-2+\dfrac{2}{\sqrt{n}}
\end{equation}

It follows that   $\langle -c_1,x_n\rangle \leqslant -1$ for $n\geqslant 4$, so that $\langle -c_1,x_n\rangle +1\leqslant 0$. This proves that $F$ is nonempty.

 Since $\langle x, Tx\rangle = \int\limits_0^1 t x^2(t)dt  \geqslant 0$ for all $x\in L_2[0,1]$, $f(x)$ is bounded from below over $F$. It is easily seen that $0\notin F$ and
 
 \begin{equation}\label{6}
 f(x)>0\quad  \text{for all}\quad x\in F.
 \end{equation}

 On the other hand,
 
 \begin{equation}
 f(x_n)=\int\limits_0^1tx_n^2(t)=\int\limits_{\dfrac{1}{n^2}}^{\dfrac{1}{n}}tndt=\dfrac{n}{2}(\dfrac{1}{n^2}-\dfrac{1}{n^4})\longrightarrow 0\quad\text{as}\quad n\to \infty.
 \end{equation}
 
 This, together with \eqref{6}, shows that the infimum of $f$ over $F$ is $0$. However, the inequality
 \eqref{6} shows that this infimum cannot be attained by any $x\in F$.
 We have shown that \eqref{QP} has no solution.

\end{example}

The following example is taken from \cite[p.~45]{BK02} to show that in the case where the set $I_1$ consists of more than one element, the condition \eqref{ConditionA} cannot be dropped from the assumptions of Theorem \ref{theo:exist1} (even in finite dimensional setting).
\begin{example}\label{ex:2} {\rm (\cite[p.~45]{BK02})} 
 Consider \eqref{eq:QP}, where ${\cal H}=\mathbb{R}^3$,
\begin{equation}\notag
\setcounter{MaxMatrixCols}{20} T=
\begin{bmatrix}
0 &0   &0 \\
0   &0 &-1  \\
0   &-1    &0 \\
\end{bmatrix}, \;\; T_1=
\begin{bmatrix}
0 &0   &0 \\
0   &1 &0  \\
0   &0    &0 \\
\end{bmatrix}, \;\;
T_2=
\begin{bmatrix}
0 &0   &0 \\
0   &0 &0  \\
0   &  0  &1 \\
\end{bmatrix}
\end{equation}
\begin{equation}\notag
c=(2,0,0),\quad c_1=(-1,0,0),\quad c_2=(-1,0,0), \quad \alpha_1=0,\quad \alpha_2=-1.
\end{equation}
We can rewrite the problem as follows:
\begin{equation*}
\begin{cases}
&\min f(x_1, x_2, x_3) := 2x_1 - 2x_2x_3  \\
& s.t. F=\{x\in \mathbb{R}^3\mid  x_2 ^2 - x_1 \leqslant 0,  x^2_3 - x_1 -1 \leqslant 0\}.
\end{cases}
\end{equation*}
 On the one hand, it is easy to check that, for this problem:
 $I_1=I=\{1,2\}$ and
 \begin{align*}
 0^+F&=\{v\in \mathbb{R}^3\mid T_1v=0, \langle c_1,v\rangle \leqslant 0, T_2v=0, \langle c_2,v\rangle \leqslant 0\}\\
 &=\{v\in \mathbb{R}^3\mid v_1\geqslant 0, v_2=0,v_3=0\},
  \end{align*}
   \begin{align*}
 \{v\in 0^+F\mid \langle v,Tv\rangle =0\}=\{v\in \mathbb{R}^3\mid v_1\geqslant 0, v_2=0,v_3=0\},
  \end{align*}
    \begin{align*}
 \{v\in 0^+F\mid \langle c_1,v\rangle =0, \langle c_2,v\rangle =0 \}&=\{v\in \mathbb{R}^3\mid v_1=v_2=v_3=0\}.
 \end{align*}
 Since the problem is setting in finite dimensional space, $\langle x,Tx\rangle $ is a Legendre form. It is easily seen that the condition \eqref{ConditionA} do not hold. On the other hand, according to \cite[p.~45]{BK02}, $f(x)$ is bounded from below over nonempty $F$ and the problem has no solution.
\end{example}
\bigskip
The following theorem shows that  for \eqref{eq:QP} with only one constraint,  condition \eqref{ConditionA} can be dropped from the assumption of Theorem \ref{theo:exist1}.
%We next  prove the solution existence  \eqref{eq:QP} with one constraint.

\begin{theorem}\label{theoL:exist3}
Consider the problem \eqref{eq:QP} where $m=1$ and $\langle x, Tx\rangle$ is a Legendre form.
Assume that the objective function is bounded from below over the nonempty feasible set.
Then \eqref{eq:QP} has an optimal solution.
\end{theorem}
{\it Proof.}
Let $f^* := \underset{x\in F}{\inf} f(x) > -\infty$. Consider the set
\begin{align*}
M&=\{v\in 0^+F\mid \langle v,Tv\rangle =0\}
=\{v\in {\cal H}\mid T_1v=0, \langle c_1,v\rangle\leqslant 0,  \langle v, Tv\rangle =0\}.
\end{align*}
 We now consider two separated cases:

If   $\langle c_1,v\rangle=0$ for all $v\in M$, then the condition \eqref{ConditionA} is satisfied. From Theorem \ref{theocomp} it follows that \eqref{eq:QP} has a solution.

Consider the case where there is  $\bar v\in M$ such that $\langle c_1, \bar v\rangle <0$.  As  $\langle c_1, \bar{v}\rangle < 0$, there exists $ t_0 \geq 0$ such that $g_1(t_0\bar{v}) < 0$.  Furthermore, the quadratic function $f$ is bounded from below over $F$.  According to Theorem \ref{S-lemma}  there exists $
\lambda \geqslant 0$ such that
\begin{align}\label{eq:theoL1}
f(x) +\lambda g_1(x)\geqslant f^* \,\; \forall x\in {\cal H}.
\end{align}

Consider the quadratic problem
\begin{align}\label{eq:theoL2}
\min \{f(x) + \lambda g_1(x)\mid \; \mbox{subject to}\;\; x\in {\cal H}\}
 \end{align}

Since $\langle x, Tx\rangle$ is a Legendre form  and $T_1$  is positive semidefinite, by Corollary 2 in  \cite[p.553]{Hestenes1}, we have  $\langle x, \left(T+\lambda T_1\right)x\rangle$ is a Legendre form. Hence,  according to Corollary \ref{cor:exist4}, the problem \eqref{eq:theoL2} has a  solution, say $x^*$.  
We have
\begin{align*}
f(x^*)+\lambda g_1(x^*)\leqslant f(x)+\lambda g_1(x)\; \; \forall x\in {\cal H}. 
\end{align*}
This gives 
$f(x^*)+\lambda g_1(x^*)\leqslant f(x)$ for all $x\in F$.
From this it follows that 
$f(x^*)+\lambda g_1(x^*)\leqslant f^*$.
Combining this with \eqref{eq:theoL1} we get
\begin{align}\label{eq:theoL3}
f(x^*)+\lambda g_1(x^*)= f^*.
\end{align}

We consider three distinguish possibilities.

 1. If  $g_1(x^*)= 0$, then $x^*\in F$ and  by \eqref{eq:theoL3},  $x^*$ is a solution of \eqref{eq:QP}.

2. Consider the case where $g_1(x^*)< 0$.   Since
\begin{align*}
 g_1(x^*-t\bar{v})& = g_1(x^*)+\frac{1}{2}t^2\langle \bar v, T_1\bar v\rangle-t\langle T_1x^*+c_1,\bar v\rangle\\
 &=g_1(x^*)- t\langle c_1, \bar{v}\rangle,
 \end{align*}
 we can choose $t^*=\frac{g_1(x^*)}{\langle c_1, \bar{v}\rangle} >0$ so that  $g_1(x^*-t^*\bar{v}) = 0$
 and then we have
  \begin{equation}\label{eq:theoL4}
x^*-t^*\bar{v}\in F.
 \end{equation}

From \eqref{eq:theoL1} we have
$f(x^*+t^*\bar{v})+\lambda g_1(x^*+ t^*\bar{v}) \geqslant f^*$
or, equivalently,
\begin{align*}
f(x^*)+t^*\langle Tx^*+c,\bar v\rangle +\lambda g_1(x^*)+\lambda t^*\langle c_1,\bar v\rangle \geqslant f^*.
\end{align*}
It follows that
$\lambda t^*\langle c_1,\bar v\rangle \geqslant -t^*\langle Tx^*+c,\bar v\rangle$.
Combining this with \eqref{eq:theoL3} and \eqref{eq:theoL4},  we have
\begin{align*}
f(x^*-t^*v)& = f(x^*)- t^*\langle T x^*+c, \bar v\rangle
 \leqslant f(x^*) +\lambda t^* \langle c_1, \bar{v}\rangle\\
&= f(x^*)+\lambda \frac{g_1(x^*)}{\langle c_1, \bar{v}\rangle} \langle c_1, \bar{v}\rangle
=f(x^*)+\lambda g_1(x^*) = f^*.
\end{align*}
It follows that $x^*-t^*v$ is a solution to \eqref{eq:QP}.

3. Consider the case where  $g_1(x^*)> 0$.   Since
\begin{align*}
 g_1(x^*+t\bar{v})& = g_1(x^*)+\frac{1}{2}t^2\langle \bar v, T_1\bar v\rangle+t\langle T_1x^*+c_1,\bar v\rangle
 =g_1(x^*)+ t\langle c_1, \bar{v}\rangle,
 \end{align*}
 we can choose $t^*=\frac{- g_1(x^*)}{\langle c_1, \bar{v}\rangle} >0$  and then we have
$g_1(x^*+t^*\bar{v})  = 0$, so 
\begin{equation}\label{eq:theoL5}
x^*+t^*\bar{v}\in F.
 \end{equation}

By \eqref{eq:theoL1}, we have 
$f(x^*-t^*\bar{v})+\lambda g_1(x^*- t^*\bar{v}) \geqslant f^*$
or, equivalently,
$$f(x^*)-t^*\langle Tx^*+c,\bar v\rangle +\lambda g_1(x^*)-\lambda t^*\langle c_1,\bar v\rangle \geqslant f^*.$$
It follows that
$-\lambda t^*\langle c_1,\bar v\rangle \geqslant t^*\langle Tx^*+c,\bar v\rangle$.
Combining this with  \eqref{eq:theoL3} and \eqref{eq:theoL5},  we have
\begin{align*}
f(x^*+t^*v)& = f(x^*)+ t^*\langle T x^*+c, \bar v\rangle
\leqslant f(x^*) -\lambda t^* \langle c_1, \bar{v}\rangle\\
&= f(x^*)+\lambda \frac{g_1(x^*)}{\langle c_1, \bar{v}\rangle} \langle c_1, \bar{v}\rangle
=f(x^*)+\lambda g_1(x^*) = f^*.
\end{align*}
It follows that $x^*+t^*v$ is a solution to \eqref{eq:QP}.
Thus we have showed that the problem \eqref{eq:QP}  has a solution. 
The proof is complete. 

\qed

In order to apply Theorem  \ref{theo:exist1}, we have to find out whether the objective function
$f(x)$ of \eqref{eq:QP} is bounded from below over $F$, or not. This is a rather difficult task.
The following theorem  gives a sufficient condition for solvability of \eqref{eq:QP}.

\begin{theorem}\label{theo2}{\rm(Eaves-type theorem)}
Consider the problem \eqref{eq:QP}, where $\langle x, Tx\rangle$ is a Legendre form.  
Suppose that  
\begin{itemize}
\item[\rm (i)] ${F}$ {\it  is nonempty;}
\item[\rm (ii)] {\it  If $v\in 0^+ F$ then $\langle v,Tv\rangle \geq 0$;}
\item [\rm (iii)] {\it  If $v\in0^+F$ and $x\in {F}$ are such that $ \langle v,Tv\rangle=0$, then $\langle Tx+c,v\rangle \geq 0$;}
\item [\rm (iv)] {\it  The condition \eqref{ConditionA} is satisfied.}
\end{itemize}
Then,  problem \eqref{eq:QP} has a solution.
\end{theorem}
{\it Proof.}
To prove that \eqref{eq:QP} has a solution under conditions (i)-(iv), by Theorem \ref{theo:exist1} it
suffices to verify that $f$ is bounded from below over $F$.

Assume that conditions (i), (ii), (iii) and (iv) are satisfied. Define $f^*=\inf\{f(x)\,:\,x\in F\}$. As $F\neq \emptyset$, we have $f^*\neq +\infty$. If $f^*>-\infty$ then the assertion of the theorem follows from the Frank-Wolfe type theorem \ref{theo:exist1}. Hence we only need to show that $f^*>-\infty$. To obtain a contradiction, suppose that $f^*=-\infty$. Then, there exists a sequence $\{y^k\}\subset F$ such that $f(y^k)\rightarrow -\infty$. There is no loss of generality in assuming that $\Vert y^k \Vert \rightarrow \infty$ as $k\to\infty$ and $f(y^k)\leqslant \frac{1}{k}$. 

Let $S_k= \{ x\in F\mid  f(x) \leqslant f(y^k) \}$. We have $y^k\in S_k$, so $S_k$ is nonempty, closed.   
 By Lemma \ref{levelset}, $S_{k}$ admits an element of minimal norm. Let $x^k\in S_{k}$ be a element of minimal norm. 
Since $f(x^k)\leqslant f(y^k)$ and $f(y^k)\rightarrow -\infty$, we have $f(x^k)\rightarrow -\infty$.
% There is no loss of generality in assuming that $\Vert x^k \Vert \rightarrow \infty$ as $k\to\infty$.
 Note that
\begin{equation}\label{eq14}
g_i(x^k)=\frac{1}{2}\langle x^k,T_ix^k\rangle +\langle c_i,x^k\rangle +\alpha_i \leqslant
0,\quad \forall i=1,..., m,\ \, \forall k\geqslant 1.
\end{equation}
Without loss of generality may assume that $\Vert x^k\Vert\not =0$ for all $k$, and 
$\frac{x^k}{\Vert x^k\Vert}\rightharpoonup \bar v$. 
By Lemma \ref{lema:rec2}, we have  $\bar v\in 0^+F$.

Since $f(x^k)\rightarrow -\infty$, we can assume that for all $k\geqslant 1$,
 \begin{equation}\label{eq15}
f(x^k)=\frac{1}{2}\langle x^k,Tx^k\rangle+\langle c,x^k\rangle\leqslant 0.
\end{equation}
Multiplying the inequality in \eqref{eq15} by $\Vert x^k\Vert^{-2}$ and letting $k\to\infty$, one has
\begin{align*}
\underset{k\to\infty}{\limsup}\frac{1}{2}\langle\Vert x^k\Vert^{-1}x^k, T\Vert x^k\Vert^{-1}x^k\rangle \leq 0.
\end{align*}
By the weakly lower semicontinuity of  $\langle x, Tx \rangle$, we have
\begin{equation}\label{eq17}
\frac{1}{2}\langle \bar v,T\bar v\rangle  \leqslant \underset{k\to\infty}{\liminf}\frac{1}{2}\langle \frac{x^k}{\Vert x^k\Vert}, T\frac{x^k}{\Vert x^k\Vert}\rangle\leqslant \underset{k\to\infty}{\limsup}\frac{1}{2}\langle\frac{x^k}{\Vert x^k\Vert}, T\frac{x^k}{\Vert x^k\Vert}\rangle \leq 0.
\end{equation}
   
From \eqref{eq17} and assumption (ii) we have 
\begin{equation}\label{eq:18}
\underset{k\to\infty}{\lim}\langle \frac{x^k}{\Vert x^k\Vert}, T\frac{x^k}{\Vert x^k\Vert}\rangle= \langle \bar v, T\bar v\rangle = 0.
\end{equation}
 Since $\langle x,Tx\rangle$ is a Legendre form, $ \frac{x^k}{\Vert x^k\Vert}\rightharpoonup \bar v$    we can deduce that $\Vert \bar v\Vert=1$.
Since $\bar v\in 0^+F$,  by \eqref{eq:18} and assumptions (iii),   we can deduce that
\begin{equation}\label{eq20}
\langle Tx^k+c,\bar v \rangle \geqslant 0. 
\end{equation}
We have shown that there exist $S_k$, $x^k$ such that the assumptions (b) and (c) in Lemma 
\ref{lema:bd} are satisfied. Since $\langle x, T x\rangle$ is a Legendre form, it is weakly lower semicontinuous. Hence,  assumption (a) in Lemma \ref{lema:bd} is satisfied. Combining this with assumption (d) we see that all the assumptions in Lemma   \ref{lema:bd} hold. Hence $\{x^k\}$ is bounded. Without loss of generality, we can assume that $x^k\rightharpoonup \bar x$ as $k\to\infty$. Since $x^k\in F$ and $F$ is weakly closed, $\bar x\in F$. By weak lower semicontinuity of $f$, $f(\bar x)\leqslant \underset{k\to\infty}{\liminf}f(x^k) = -\infty$,
   which is impossible. The proof is complete. 
   
\qed

The following example shows that the condition \eqref{ConditionA} from assumption of Theorem \ref{theo2}  is  not necessary for the solution existence of \eqref{eq:QP}.

\begin{example}\label{exam_5}
 Consider the  programming problem
\begin{align}\label{prob_exam_5}
\begin{cases}
&\mbox{min } f(x) =\frac{1}{2} \langle x, T x\rangle +\langle c, x\rangle, \\
&\mbox{subject to } x\in \ell_2: \frac{1}{2}\langle x, T_i x\rangle +\langle c_i, x\rangle+\alpha_i\leqslant 0 ,\;\; i=1,2.
\end{cases}
\end{align}
where $Tx =(0,-x_2,x_3,x_4,\ldots)$, $c =(0, 1,0,0,\ldots)$, $T_1x =(0, x_2, x_3, \ldots)$, $c_1 = (0,0,0,\ldots)$, $\alpha_1 = 0$, $T_2x=(0, 0, x_3, 0, 0 \ldots)$, $c_2=(1, 0, 0 \ldots)$, $\alpha_2 = 0$.

 It is clear that $\langle x, T x\rangle$ is a Legendre form and 
\begin{align*}
&0^+F = F =\{v= (v_1,0,0,\ldots)\in \ell_2\mid v_1\leqslant 0\},\\ 
&\{v\in 0^+F\mid \langle v, Tv\rangle = 0\}  = \{v= (v_1,0,0,\ldots)\in \ell_2\mid v_1\leqslant 0\}.
\end{align*}
Since the $f(x) = 0$ on $F$, the solution set of \eqref{prob_exam_5} coincides with $F$. 
For $\bar v := (-1,0,...,0,\ldots)\in \{v\in 0^+F\mid \langle v, Tv\rangle = 0\} $ we have   $\langle c_2,\bar v\rangle \neq 0$. Hence the condition \eqref{ConditionA} does not hold for \eqref{prob_exam_5}.
\end{example}

\begin{remark} 
Note that either $T_i = 0$ for all $i=1,\ldots, m$ or $c_i = 0$ for all $i\in I_1$,  then  condition \eqref{ConditionA}is automatically statisfied.  In each of the cases, it is to show that assumptions   $\rm (ii)$, $\rm (iii)$ are a necessary and
sufficient  for the solution existence of \eqref{eq:QP}, provided that $F\neq \emptyset$. 
\end{remark}
%\subsection{Existence theorem involving  compact operators}

%In the remainder of this section we will give some existence results by using the notion of compact linear operator.
In the remainder of this section we give an existence result of the solution
for \eqref{eq:QP} under the assumption that all the operators corresponding to quadratic forms are compact operators  with closed range. This  means that all the operators corresponding to quadratic forms in \eqref{eq:QP}  have finite dimensional ranges (see, \cite[Theorem 4.18]{Rudin}). Note that this assumption is very restrictive but by using this assumption we can investigate
  the solution existence for a class of  \eqref{eq:QP} problems where the quadratic form in objective function is not a Legendre form. The next statement may be seen as a complement to Theorem \ref{theo:exist1}.

\begin{theorem}\label{theocomp}{\rm (Frank-Wolfe-type theorem 2)}
Consider the problem \eqref{eq:QP} and suppose that
\begin{itemize}
\item [\rm (i)]  $T$ and  $T_i$ ($ i = 1,2,\dots, m$) are compact with closed range;
\item [\rm (ii)]  the objective function $f$ is bounded from below over the nonempty $F$;
\item [\rm (iii)] the condition \eqref{ConditionA} is satisfied.
\end{itemize}
Then, \eqref{eq:QP} has a solution.
\end{theorem}
{\it Proof.}  
Let $f^*=\inf \{f(x)\mid x\in F\}$. For each positive integer $k$, put $S_k=\{x\in F \mid f(x)\leqslant f^*+\frac{1}{k}\}$.
By assumption $f^* >-\infty$, $S_k$ is nonempty and closed. 
By Lemma \ref{levelset}, $S_k$ admits an element of minimal norm. Let $x^k\in S_k$ be a element of minimal norm. Then, we have
\begin{align}
\label{eq:theocomp1}\; f(x_k)&=\frac{1}{2}\langle  x^k, Tx^k\rangle+\langle c, x^k\rangle   \leqslant f^*+\frac{1}{k},\\
 \label{eq:theocomp2} \; g_i (x_k)&=\frac{1}{2}\langle  x^k, T_i x^k\rangle+\langle c, x^k\rangle +\alpha_i  \leqslant 0\; \; i=1,\dots, m.
\end{align}
  Consider sequence $\{x^k\}$. We now prove that $\{x^k\}$ is bounded. On the contrary, suppose that $\{x^k\}$ is unbounded. Without loss of generality we may assume that
$\Vert x_k\Vert\to\infty$ as $k\to\infty$, $\Vert x_k\Vert \neq 0$ for all $k$.
Put $v^k:= \frac{x^k}{\Vert x^k\Vert}$, one has $\Vert v^k\Vert = 1$,
then there exists a subsequence of $v^k$ weakly converges  to $\bar v$.
Without loss of generality we can assume that $v^k \rightharpoonup \bar v$ as $k\to\infty$. Since $T$ is compact, by Hilbert theorem (see, \cite[p. 261]{Ioffe}), $\langle x,Tx\rangle $ is weakly continuous . By an argument analogous  to those in the proof of Theorem \ref{theo:exist1} we can deduce
\begin{align*}
&\bar v\in 0^+F,\, 
\langle \bar v, T\bar v\rangle  = 0.
\end{align*}

Consider two cases:

Case 1:  $\bar v\neq 0$. In this case, from the assumptions of the Theorem and from the above, the assumptions of  Lemma \ref{lema:bd} are satisfied. By Lemma \ref{lema:bd}, $\{x^k\}$ is bounded.

Case 2:  $\bar v= 0$.  Then,
\begin{align*}
 T\bar v = 0,\;\;  \langle c, \bar v\rangle = 0,\;\; T_i\bar v = 0,\; \; \langle c_i, \bar v \rangle =0,\quad \forall i=1,2,\cdots, m.
\end{align*}

Let ${\cal L} =\underbrace{{\cal H}\oplus \dots \oplus {\cal H}}_{(m+1)\; times}\oplus \mathbb{R}^{m+1}$, where $\oplus$ denotes the direct sum of Hilbert spaces, and let us write $\langle . , . \rangle _{\cal L}$ and $\Vert .\Vert_{\cal L}$ for the scalar product and the norm on ${\cal L}$, respectively.

Let $A: { \cal H}\to {\cal L}$ be defined by
$$Ax=(Tx, T_1x,\dots , T_m x, \langle c,x\rangle, \langle c_1,x\rangle, \dots, \langle c_m,x\rangle).$$
 Since $T$ and  $T_i$ ($i=1,\dots, m$) are continuous with closed range,  so is $A$. For each $k$, consider the linear system
\begin{align}\label{eq:theocomp3}
Ax=Ax^k.
\end{align}
Since $A$ is continuous with closed range, there exist the continuous pseudoinverse $A^+$ of $A$ and a solution $\bar x^k$ to \eqref{eq:theocomp3} such that $\bar x^k=A^+Ax^k$ (see, for instance, \cite{Luenberger}, p. 163). Therefore, there exists $\rho > 0$ , depending on $A$, such that
$\Vert \bar x^k\Vert\leqslant \rho\Big( \Vert Ax^k\Vert_{\cal L} \Big)$.
 This gives
\begin{align*}
\Vert \bar x^k\Vert\leqslant \rho\Big( \Vert Tx^k\Vert +\sum_{i=1}^m \Vert T_ix^k \Vert+ \vert \langle c,x^k\rangle \vert+ \sum_{i=1}^m \vert \langle c_i, x^k\rangle \vert\Big).
\end{align*}

 By \eqref{eq:theocomp3}, $A\bar x^k=Ax^k$, we can check that $\bar x^k\in S_k $ and
 \begin{equation*}
f(\bar x^k)= f(x^k)\leqslant f^*+\frac{1}{k}
\end{equation*}
Since $x^k$ is the  minimal norm element in $S_k$, we have
 \begin{equation*}
\Vert x^k\Vert\leqslant \Vert \bar x^k\Vert\leqslant \rho \Big( \Vert Tx^k\Vert +\sum_{i=1}^m \Vert T_ix^k \Vert+ \vert \langle c,x^k\rangle \vert+ \sum_{i=1}^m \vert \langle c_i,x^k\rangle \vert\Big)\;\; \forall k.
\end{equation*}
Dividing both sides of the above inequality by $\Vert x^k\Vert$, letting $k\to \infty$ and by the compactness of $T$, one has
\begin{align*}
1\leqslant\rho\Big( \Vert T\bar v\Vert +\sum_{i=1}^m \Vert T_i\bar v \Vert+ \vert \langle c,\bar v\rangle \vert+ \sum_{i=1}^m \vert \langle c_i,\bar v\rangle \vert\Big).
\end{align*}
This contradicts the fact that  $T \bar v = 0$, $\langle c, \bar v\rangle=0$, $T_i\bar v=0$, $c_i\bar v=0$, $i=i,\dots , m$. 

Thus $\{x^k\}$ is bounded. Therefore,  $\{x_k\}$ has a weakly convergent subsequence. Without loss of generality,
we may assume that $x^k\rightharpoonup \bar x$ as $k\to\infty$. Since $F$ is weakly closed, $x^k\in F$, we have $\bar x \in F$. Since $ T$ is compact, by Theorem Hilbert (see \cite[p.261]{Ioffe}), $\langle x, Tx\rangle$  is  weakly continuous. Hence,  by  \eqref{eq:theocomp1},
\begin{align*}
f(\bar x)=\frac{1}{2}\langle \bar x, T\bar x\rangle+\langle c, \bar x\rangle&  = \underset{k\to\infty}{\lim}\Big(\frac{1}{2}\langle x^k, T x^k\rangle+ \langle c, x^k\rangle\Big)\\
&\leqslant \underset{k\to\infty}{\lim}(f^*+\frac{1}{k})
= f^*.
\end{align*}
From the above it follows that  $\bar x$ is a solution of \eqref{eq:QP}. This proof of theorem is complete.

\qed

 Note that in the following consequences of Theorem \ref{theocomp}  the condition \eqref{ConditionA} is automatically satisfied.

\begin{corollary}\label{cor :exist1} Consider the quadratic programming problem under linear constraints ${\rm (QPL)}$ {\rm(}i.e. \eqref{eq:QP} with $T_i=0$ for all $i=1,...,m${\rm)}, where $T$ is compact operator with closed range. Suppose that $f(x)$ is bounded from below over nonempty $F$.
 Then,  problem ${\rm (QPL)}$  has a solution.
\end{corollary}
\begin{corollary}\label{cor :exist1} Consider the  linear
programming problem  ${\rm (LP)}$ {\rm(}i.e. \eqref{eq:QP} with $T=0,\; T_i=0$ for all $i=1,...,m${\rm)}. Suppose that $f(x)$ is bounded from below over nonempty $F$.
 Then,  problem ${\rm (LP)}$  has a solution.
\end{corollary}
\begin{corollary}\label{cor :exist2} Consider the problem \eqref{eq:QP}, where $T$ and $T_i$ for all $i=1,...,m$, are compact operators with closed range. Suppose that  $c_i=0$ for all $i\in I_1$ and $f(x)$ is bounded from below over nonempty $F$.
 Then,  problem \eqref{eq:QP} has a solution.
\end{corollary}

\begin{corollary}\label{cor :exist3} Consider the problem \eqref{eq:QP}, where $T$  and $T_i$ for all $i=1,...,m$, are compact operators with closed range. Suppose that $\{v\in 0^+F\mid \langle v,Tv\rangle =0\}\subset \{0\}$  and $f(x)$ is bounded from below over nonempty $F$.
 Then,  problem \eqref{eq:QP} has a solution.
\end{corollary}

\begin{corollary}\label{cor :exist4}
Let $T$ is compact operator with closed range. Suppose that quadratic function  $f(x)=\frac{1}{2}\langle x, Tx\rangle+\langle c, x\rangle$ is bounded from below over Hilbert space ${\cal H}$. Then, there exists an $x^*\in {\cal H}$ such that $f(x^*)\leqslant f(x)$ for all $x\in {\cal H}$.
\end{corollary}

\begin{remark}
{\rm If ${\cal H}$ is of finite dimension then, any continuous operator $T$ on ${\cal H}$  is compact with closed range and $\langle x, Tx\rangle$ is a Legendre form. Therefore,  in the finite dimensional setting Theorem \ref{theo:exist1} and Theorem \ref{theocomp} are identical.}
\end{remark}

\section{Conclusions}
In this paper we consider quadratic programming problems in Hilbert spaces of  arbitrary dimension and propose conditions for the solution existence of quadratic programming problems whose constraint set is defined by finitely many convex quadratic inequalities. Our results extend some previous existence results for  nonconvex quadratic programming problems in finite dimensional setting, including the Frank-Wolfe theorem for linearly constrained quadratic programming, to the inﬁnite dimensional setting.

In finite dimensional setting,  Luo and  Zhang \cite[ Theorem 2]{LuoZhang99} showed that quadratic programming problem \eqref{eq:QP}  whose the objective function is bounded from below over the nonempty constraint set defined  by finitely many convex quadratic inequalities, in which at most one is nonlinear, always has a solution. From this it follows that in the case where $I_1$ consists of at most one element then, the condition \eqref{ConditionA} can be dropped from the assumptions of Theorem \ref{theo:exist1}. Note that if $I_1$ consists of more than one element then the condition  \eqref{ConditionA} can not be dropped (see example \ref{ex:2}).

In connection with Theorem \ref{theo:exist1},  Theorem \ref{theoL:exist3} and Theorem \ref{theocomp}, one may ask whether in the case where  $I_1$ is a singleton  and $I_0$ is nonempty, the condition \eqref{ConditionA} can be dropped from the assumptions of Theorem \ref{theo:exist1} and Theorem \ref{theocomp}?. This is equivalent to the question: Can we prove the same result as in \cite[ Theorem 2]{LuoZhang99} for a quadratic programming problem \eqref{eq:QP} in Hilbert spaces of  arbitrary dimension?
\section*{Acknowledgments}
This research is funded by Vietnam National Foundation for Science and Technology Development (NAFOSTED) under grant number 101.01-2014.39. 
The authors would like to thank Prof. N. D. Yen for valuable remarks and suggestions.

\end{document}